\documentclass[12pt]{article}
\usepackage{amssymb}
\def\bb{{\mathfrak b}}
\def\dd{{\mathfrak d}}
\def\cc{{\mathfrak c}}
\def\pp{{\mathfrak p}}

\def\proof{\par\noindent Proof\par\noindent}
\def\poset{{\mathbb P}}
\def\forces{{| \kern -2pt \vdash}}
\def\force{\forces}
\def\res{\upharpoonright}
\def\qed{\par\noindent QED\par}
\def\daglamb{$(\dagger)$-$\lambda'$-set }

\def\rmand{\mbox{ and }}
\def\gameone{{G_{{\cal O,P}}(X_F,\infty)}}
\def\gamefin{{G_{{\cal O,P}}^f(X_F,\infty)}}
\def\gameonegen{{G_{{\cal O,P}}(X,x)}}
\def\gamefingen{{G_{{\cal O,P}}^f(X,x)}}
\def\gamegam{{G_{{\cal F,C}}^{\gamma}(X)}}
\def\cantorxplus{X_F}
\newtheorem{theorem}{Theorem}[section]
\newtheorem{lemma}[theorem]{Lemma}

\newtheorem{question}[theorem]{Question}

\begin{document}

\begin{center}
{\large On $\lambda'$-sets}
\end{center}

\begin{flushright}
Arnold W. Miller\footnote{
Thanks to the Fields Institute for Research in Mathematical Sciences
at the University of Toronto for their support during the time this paper
was written and to Juris Steprans who directed the special program in set
theory and analysis.
}
\end{flushright}

\begin{center}
Abstract
\end{center}

\begin{quote}
A set $X\subseteq 2^\omega$ is a $\lambda'$-set iff for every countable set
$Y\subseteq 2^\omega$ there exists a $G_\delta$ set $G$ such that 
$(X\cup Y)\cap G=Y$.  In this paper we prove two forcing 
results about $\lambda'$-sets.
First we show that it is consistent that every $\lambda'$-set is
a $\gamma$-set.  Secondly we show that is independent whether or
not every \daglamb is a $\lambda'$-set.  
\end{quote}

\section{$\lambda'$-sets and $\gamma$-sets}


A set $X\subseteq 2^\omega$ is a $\lambda'$-set iff for all countable
$A\subseteq 2^\omega$ there exists a $G_\delta$ set $G$ such that
$$(X\cup A)\cap G=A$$ 
An $\omega$-cover of $X$ is a countable set of open sets such that
every finite subset of $X$ is contained an element of the cover.
A $\gamma$-cover of $X$ is a countable sequence of open subsets of 
$X$ such that every element of $X$ is in all but countably many 
elements of the sequence. 

\medskip\noindent
Define. $X$ is a $\gamma$-set iff any $\omega$-cover of $X$ contains a
$\gamma$-cover of $X$.

\medskip In this section we answer a question of Gary Gruenhage 
who asked if there is always a $\lambda'$-set which is not a $\gamma$-set. We
answer this in the negative.

It is well known (see Gerlitz and Nagy \cite{GN}) that MA($\sigma$-centered)
implies that every set of reals of cardinality less than the continuum is a
$\gamma$-set.  The standard model for MA($\sigma$-centered) (see Kunen and Tall
\cite{KT}) is obtained as follows:

Suppose that $M$ is a countable standard model of  ZFC+CH and we iterate
$\sigma$-centered forcings of size $\omega_1$ in $M$ with a finite support
iteration of length $\omega_2$.  In the final model $M_{\omega_2}$,  we have
that MA($\sigma$-centered) is true and the continuum is $\omega_2$.

\begin{theorem} 
In the standard model for MA($\sigma$-centered) every 
$\lambda'$ set has cardinality $\leq \omega_1$.
and (it follows
from MA($\sigma$-centered)) every set of size $\omega_1$ is
a $\gamma$-set. Hence, in this model, every $\lambda'$-set is
a $\gamma$-set.
\end{theorem}
\proof

We will use the following Lemma in our proof.

\begin{lemma}
 Suppose that $\poset$ is a $\sigma$-centered forcing
such that 
 $$\forces \tau\in {2^\omega}$$  
Then there
exists a countable set $A\subseteq 2^\omega$ in the ground model
such that for every $p\in\poset$ and open set $U\supseteq A$
coded in the ground model
there exists $q\leq p$ such that $q\forces \tau\in U$.
\end{lemma}
\proof

To prove the Lemma we will use the following Claim.

\bigskip
Claim.  Suppose $\Sigma\subseteq \poset$ is a centered subset. Then
there exists $x\in 2^\omega$ such that for every $p\in\Sigma$ and
for every $n<\omega$ there exists $q\leq p$ such that
 $$p\forces \check{x}\res n= \tau\res n.$$

\noindent pf: Otherwise by the compactness of $2^\omega$ there exists a finite
set  $$\{p_m:m<N\}\subseteq \Sigma \rmand \{s_m:m<N\} \subseteq 2^{<\omega}$$
such that $\{[s_m]:m<N\}$ covers $2^\omega$ and for each $m<N$ we have that
       $$p_m\forces \tau\notin [s_m].$$ 
But this is a contradiction since there exists
some $p\in\poset$ below all of the $p_m$. This proves the Claim.

\bigskip

Let $\poset=\bigcup_{n<\omega} \Sigma_n$ be a sequence
of centered sets.  Then for each $n$ there exists $x_n\in 2^\omega$
such that for every $p\in \Sigma_n$ and for every $m\in \omega$
there exists $q\leq p$ such that 
 $$ q\forces \check{x}_n\res m=\tau\res m.$$
Now let $A=\{x_n:n<\omega\}$.  This proves the Lemma.
\qed

Suppose $X\subseteq 2^\omega$ is a $\lambda'$-set  in $M_{\omega_2}$.
For each $\alpha\leq\omega_2$ define 
 $$X_\alpha=X\cap M_\alpha$$
By a standard Lowenheim-Skolem argument we can find $\alpha<\omega_2$
such that
\begin{enumerate} 
  \item $X_\alpha\in M_\alpha$ and
  \item for every countable $A\subseteq 2^\omega$ which is
    in $M_\alpha$ there exists a $G_\delta$-set $G$ coded in 
    $M_\alpha$ such that
           $$(X_{\omega_2}\cup A)\cap G=A$$ 
\end{enumerate}
We claim that $X=X_{\omega_2}=X_\alpha$ and hence has cardinality
$\leq\omega_1$.  Suppose that $\tau$ is any term for an element of $2^\omega$
in $M_{\omega_2}$.  Since $\tau$ is added at some latter stage $\beta$ with
$\alpha\leq \beta <\omega_2$  and the iteration of $\sigma$-centered forcings
of length $<\omega_2$ is $\sigma$-centered, it  follows that $\tau$ is added by
a $\sigma$-centered forcing over $M_\alpha$.  Let $A\subseteq2^\omega$ be the
countable set given by the Lemma.  By the Lemma it follows that $\tau$ must be
an element of any $G_\delta$ set coded in $M_\alpha$ which contains $A$.  
Using item (2) above we see that $\tau$ must be in $A$ if it is in
$X_{\omega_2}$. Therefore $X_{\omega_2}\setminus X_\alpha=\emptyset$.

\qed
 
\bigskip
Remark.  This argument is similar to the proof that there are no
$\lambda'$-sets of size $\omega_2$ in Laver's model, see Miller \cite{survey2}.

\bigskip

Remark.  A set of reals $X$ is a $\lambda$-set iff every countable subset of
$X$ is a relative $G_\delta$.   In ZFC we must always have a  $\lambda$-set
which is not a $\gamma$-set. 
To see this let 
 $$X=\{f_\alpha\in \omega^\omega: \alpha<\bb\}$$ 
be well-ordered by eventual dominance and unbounded. Then  Rothberger
\cite{rothberger} (or see Miller \cite{survey1}) showed that $X$ is a
$\lambda$-set. However $X$ is not a $\gamma$-set as is witnessed by the
sequences of $\omega$-covers 
 $${\cal U}_m=
 \{U_n^m:n\in\omega\}\mbox{ where } U_n^m=\{f\in\omega^\omega: f(m)<n\}.$$ 

In fact the set $X$ is a $\lambda'$-set with respect to
$\omega^\omega$. This follows from the following lemma.

\begin{lemma}(Rothberger)
Suppose $Z_\beta=\{f_\alpha:\alpha<\beta\}\subseteq \omega^\omega$ is
well-ordered by eventual dominance,
and $A\subseteq \omega^\omega$ is countable and for
every $g\in A$ there exists $\alpha<\beta$ such that
$\exists^\infty n\; g(n)<f_\alpha(n)$.  Then 
there exists a $G_\delta$ set $G$ with
 $$G\cap (Z_\beta\cup A)=A$$  
\end{lemma}
\proof
This is proved by induction on $\beta$. 
and assume the lemma is true for all $\delta<\beta$.  If
$\beta$ is a successor ordinal, then the induction is trivial.  

\noindent Case 1. $\beta$ is a limit ordinal of uncountable cofinality.

Find $\delta_0<\beta$ so that for each $g\in A$
$\exists^\infty n\;\; g(n)<f_{\delta_0}(n)$.
Then by induction there exists a $G_\delta$ set $G$ with
$$G\cap (Z_{\delta_0}\cup A)=A$$
Let $H=\{g\in\omega: \exists^\infty n\;\; g(n)<f_{\delta_0}(n)\}$
Then $H$ is a $G_\delta$ set containing $A$ and missing
$Z_\beta\setminus Z_\delta$ and so
$$(G\cap H)\cap (Z_\beta\cup A)=A$$

\noindent Case 2. $\beta$ is a limit ordinal of countable cofinality.

Let $\beta_n$ be an increasing $\omega$-sequence with limit $\beta$
and let 
$$A_n=\{g\in A: \exists^\infty m\;\; g(m)<f_{\beta_n}(m)\}$$
By inductive assumption there exists $G_\delta$ sets 
$G_n$ so that
$$G_n\cap (Z_{\beta_n}\cup A_n)=A_n$$
Define
$$G_n^*=G_n\cup\{g\in\omega^\omega:\exists^\infty m\;\; f_{\beta_n}(m)
\leq g(m)\}$$
Note that $G_n^*$ is a $G_\delta$ set which contains $A$ but still
$$G_n^*\cap (Z_{\beta_n}\cup A_n)=A_n$$
Define $G=\cap_{n<\omega}G_n^*$.  Then $G$ is a $G_\delta$-set with
$$G\cap (Z_\beta\cup A)=A$$
\qed

\bigskip

Remark. A Hausdorff gap is an example of a $\lambda'$ set of cardinality
$\omega_1$.  $\gamma$-sets have strong measure zero and Laver \cite{laver}
proved that it consistent that every strong measure zero set is countable.

Suppose there exists $X,Y\subseteq 2^\omega$ such that $|X|=|Y|$ and $X$ is a
$\lambda'$-set and $Y$ is not a $\gamma$-set.  Then there exists $Z$ which is
a  $\lambda'$-set and not a $\gamma$-set. To see this let
$X=\{x_\alpha:\alpha<\kappa\}$ and  $Y=\{y_\alpha:\alpha<\kappa\}$.   Put 
$Z=\{(x_\alpha,y_\alpha):\alpha<\kappa\}$. The first $\kappa$ for which 
MA($\sigma$-centered) fails is $\pp$  (Bell \cite{bell}) and $\pp$ is
also the size of the smallest non $\gamma$-set.  Hence any model where every
$\lambda'$-set is $\gamma$-set and $\cc\leq\omega_2$
must satisfy MA($\sigma$-centered) and $\cc=\omega_2$.

\bigskip

Remark. Gruenhage and Szeptychi \cite{gruen} were 
interested in obtaining a set of reals
$X\subseteq 2^\omega$ which is $\gamma$-set and not a $\lambda'$-set  because
of the following two topological games.

\bigskip
Let $X$ be a topological space and $x\in X$.

\bigskip
Game: $\gameonegen$: On round $n$ player ${\cal O}$ chooses an open
neighborhood $U_n$ of $x$ and player ${\cal P}$ chooses a point
$p_n\in U_n\setminus\{x\}$. Player ${\cal O}$ wins iff the sequence
$p_n$ converges to $x$.

\bigskip
Game: $\gamefingen$:  The same except
we allow player ${\cal P}$ to choose a finite set of points
$P_n\subseteq U_n\setminus\{x\}$ on his move and ${\cal O}$ 
wins iff $\cup_{n<\omega}P_n$ converges to $x$.

\bigskip
It is not hard to check that player ${\cal O}$ has a winning strategy in
$\gameonegen$ iff player ${\cal O}$ has a winning strategy in $\gamefingen$.
Also if player ${\cal P}$ has a winning strategy in $\gameonegen$, then it is
a winning strategy in $\gamefingen$.

Given $X\subseteq 2^\omega$ consider the topology on 
$2^{<\omega}\cup {\infty}$
generated by 
 \begin{enumerate}
   \item $\{\sigma\}$ for each  $\sigma\in 2^{<\omega}$ and
   \item $\{\infty\}\cup (2^{<\omega}\setminus \cup\{x\res n: n<\omega\})$
     for each $x\in X$.
 \end{enumerate} 
Let $\cantorxplus$ denote this countable topological space.

\bigskip\noindent Gruenhage \cite{gruen1}, Nyikos \cite{nyikos},
Sharma \cite{sharma}, and Gruenhage and Szeptycki \cite{gruen}
can be combined to show that:
 
$X$ is not a $\gamma$-set iff
player ${\cal P}$ has a winning strategy in $\gamefin$.

If $X$ is a $\lambda'$-set, then
${\cal P}$ has no winning strategy in $\gameone$. 

\bigskip\noindent  Hence, if there is a set $X$ which is a $\lambda'$-set and
not a $\gamma$-set, then ${\cal P}$ has a winning strategy in $\gamefin$ but
not in $\gameone$.

\bigskip\noindent Dow \cite{dow} results imply that in Laver's model
\cite{laver}:

$X$ is a $\lambda'$-set iff
${\cal P}$ has no winning strategy in $\gameone$.

\bigskip\noindent But, it also consistent that they are not the same. In Galvin
and Miller \cite{galvin} it is shown that assuming MA($\sigma$-centered) there
is a $\gamma$-set $X$ which is concentrated on a countable subset of itself. 
Hence ${\cal P}$ has no winning strategy in $\gamefin$ hence none in
$\gameone$, but $X$ is not a $\lambda'$-set.

\begin{question}
Is it consistent with ZFC that for every $X\subseteq 2^\omega$
that

${\cal P}$ has no winning strategy in $\gameone$

iff

${\cal P}$ has no winning strategy in $\gamefin$?

\end{question}

\bigskip
To better see the connection with $\gamma$-sets consider the following
game:

\bigskip
Game: $\gamegam$: Two players ${\cal F}$ finite and ${\cal C}$ clopen
alternate plays as follows. On round $n$ player ${\cal F}$ plays
a finite set $F_n\subseteq X$ and player ${\cal C}$ responds with a clopen
set $C_n$ in $2^\omega$ with $F_n\subseteq C_n$.  Player ${\cal F}$ wins iff 
$\langle C_n:n<\omega\rangle$ is a $\gamma$-cover of $X$, ie. for
all $x\in X$ for all but finitely many $n$ we have $x\in C_n$.

\bigskip

This game is exactly the same as $\gamefin$. A neighborhood basis for
$\infty$ in $\cantorxplus$ consists of sets of the form
$2^{<\omega}\setminus\{ x\res n: x\in F, n<\omega\}$ for $F\subseteq X$
finite.  So we can regard ${\cal O}$ as player ${\cal F}$ playing a finite
subset of $X$. Instead of ${\cal P}$ playing a finite set
$P_n\subseteq 2^{<\omega}$  just regard him as ${\cal C}$ playing the clopen
set $$C_n=2^\omega\setminus \bigcup\{[s]:s\in P_n\}.$$

\begin{theorem}(Gruenhage, Szeptycki, Nyikos)
For $X\subseteq 2^\omega$ the following are equivalent:
\begin{enumerate}
 \item  $X$ is not a $\gamma$-set
 \item  ${\cal C}$ has a winning strategy in $\gamegam$.
\end{enumerate}
\end{theorem}
\proof
Suppose $X$ is is not a $\gamma$-set and let ${\cal U}$ be an $\omega$-cover
with no $\gamma$-subcover. Without loss of generality we may assume the
elements of ${\cal U}$ are clopen.  Given any $F_n$ let ${\cal C}$ choose
$C_n\in{\cal U}$ with $F_n\subseteq C_n$.
Then since $\langle C_n:n<\omega\rangle$ is not a $\gamma$-cover,
${\cal C}$ wins.

For the other direction suppose Player ${\cal C}$ has a winning strategy 
$\tau$ in $\gamegam$.  Construct
$\langle F_s,C_s: s\in \omega^{<\omega}\rangle$ so that
\begin{enumerate}
 \item for each $s\in \omega^{<\omega}$ the set
    ${\cal U}_s=\{C_{sn}:n<\omega\}$ is an $\omega$-cover of $X$ and
 \item for each $s\in \omega^{<\omega}$ and the set
    $C_s$ is the response of player ${\cal C}$ using the strategy $\tau$
    against the play $F_{s\res 1},F_{s\res 2},\ldots,F_s$.
\end{enumerate}
To do this just let
  $${\cal U}_s=\{C: \exists F \;\;
C=\tau(F_{s\res 1},F_{s\res 2},\ldots, F_{s},F)\}$$
This is countable since there are only countably many clopen sets
and by the rules of the game it must be an $\omega$-cover.  
For each element of ${\cal U}_s$ choose a witness $F$.  

Suppose for contradiction that $X$ is a
$\gamma$-set. It is well known (Gerlits and Nagy \cite{GN})
that for a $\gamma$ set
$X$ that given a sequence of
$\omega$-covers, we may choose one element of
each to get a $\gamma$-cover.  This is denoted
$X\in S_1(\Omega,\Gamma)$.  Hence we may choose $C_{sn_s}$ for
each $s\in \omega^{<\omega}$ such that every $x\in X$ is in
all but finitely many $C_{sn_s}$.  But now just look at the branch
$$m_0, m_1,m_2,\ldots \mbox{ where } m_0=n_{\langle\rangle},\ldots,
m_{k+1}=n_{\langle m_0,m_1,m_2,\ldots,m_k\rangle}$$
But
 $$F_{\langle m_0\rangle}, C_{\langle m_0\rangle},\ldots, 
 F_{\langle m_0,m_1,\ldots, m_k\rangle},
 C_{\langle m_0,m_1,\ldots, m_k\rangle},\ldots $$
is a play using the strategy $\tau$ with yields a $\gamma$ cover.
This is a contradiction.

\qed

\section{\daglamb}


In this section we answer Problem 2.12 from Nowik and Weiss 
\cite{NW} which asks basically whether it is true that
every \daglamb is a $\lambda'$-set.

Definition. For any $f\in\omega^\omega$
 $$G_f=\{a\in[\omega]^\omega\subseteq 
2^\omega:  
\forall n \exists m>n \;\;a_n < f(n)\}\;\;\;
a=\{a_0<a_1<\cdots\}$$

Definition.  A set $X\subseteq 2^\omega$ is
a \daglamb iff 
for every $f\in\omega^\omega$ we have $X\cap G_f$
is a $\lambda'$-set.

\begin{theorem}\label{dagger1}
Suppose that the continuum hypothesis is true or even just
$\bb=\dd$. Then there exists a 
\daglamb   which is not a $\lambda'$-set.
\end{theorem}

\begin{theorem}\label{dagger2}
In the Cohen real model (Cohen's original model for not CH) every
\daglamb   is a $\lambda'$-set.
\end{theorem}

\bigskip
\noindent Proof of Theorem \ref{dagger1}

Assume CH.  Let $\{f_\alpha\in\omega^\omega:\alpha<\omega_1\}$
be a scale.  That is, for $\alpha<\beta$ we have that
$f_\alpha<^* f_\beta$ and for all $g\in\omega^\omega$ there exists
$\alpha<\omega_1$ such that $g<^* f_\alpha$. 
We may also assume that the $f_\alpha$ are strictly increasing.
Let $X\subseteq [\omega]^\omega$
be the set of ranges of the elements of the scale. Then for
any $g\in\omega^\omega$ we have that $G_g\cap X$ is countable and hence
a $\lambda'$-set. On the other hand $X$ is not a $\lambda'$-set because
of the countable set $[\omega]^{<\omega}$. If $U\subseteq P(\omega)$ is
an open set containing $[\omega]^{<\omega}$, then 
$P(\omega)\setminus U$ is a compact subset of $[\omega]^{\omega}$.
If we identify $\omega^\omega$ with $[\omega]^\omega$
this means that there exists $f\in\omega^\omega$ such that
for all $g\in K$ we have $\forall n\;\; g(n)<f(n)$.
It follows that for all but
countably many $\alpha$ we have that the range($f_\alpha)\in U$.

The proof using $\bb=\dd$ is similar.  Start with a scale indexed
by $\bb$ and note that any set $Y\subseteq P(\omega)$ of size less
than $\bb$ is a $\lambda'$-set (this is due to Rothberger, see
the proof of Lemma \ref{lemdagg2}).  

\qed

\bigskip
\noindent Proof of Theorem \ref{dagger2}

Assume that $M$ is a countable transitive standard model of ZFC+CH. 

For any $\alpha\leq\omega_2^M$ let $\poset_\alpha$ be the finite 
partial functions from $\alpha$ into 2.
We claim that for any $G$ a
$\poset_{\omega_2}$-generic filter over $M$ that in the model $M[G]$ every
\daglamb is a $\lambda'$-set.  

\begin{lemma}\label{lemdagg1}
Suppose $N$ is a countable standard model of ZFC+CH, $\poset$ is a countable
poset in $N$, and 
 $$N\models X\subseteq \omega^\omega \mbox{ is unbounded in } \leq^*$$
Then for any $G$ which is $\poset$-generic over $N$ we have that
 $$N[G]\models X \mbox{ is unbounded in } \leq^*$$
\end{lemma}
\proof
Let $\{g_\alpha:\alpha<\omega_1^N\}$ be a scale in $N$. 
Working in $N$ choose $f_\alpha\in X$
so that 
 $$\exists^\infty n \;\;f_\alpha(n)>g_\alpha(n)$$  
Note that for
every $g\in\omega^\omega\cap N$ there exists $\alpha<\omega_1$ such that
$$\forall\beta>\alpha\;\; \exists^\infty n \;\;f_\beta(n)>g(n).$$

Suppose for contradiction that for some $g\in N[G]\cap\omega^\omega$
and all $\alpha<\omega_1$ we have that $g\geq^* f_\alpha$.  Then for
some $\Sigma\in[\omega_1]^{\omega_1}$ and $n<\omega$ we have
that 
 $$\forall m>n\;\;\forall \alpha\in\Sigma\;\; f_\alpha(m)\leq g(m)$$
Let $q\in G$ force this fact.
Now since $\poset$ is a countable poset, there exists some $p\in G$
with $p\leq q$ such that 
 $$\Gamma=\{\alpha<\omega_1: p\forces \alpha\in\dot{\Sigma}\}$$
is uncountable (and by definability of forcing it is in $N$).   But note that
$\{f_\alpha:\alpha\in\Gamma\}$ is unbounded and so for some $m>n$
the set $\{f_\alpha(m):\alpha\in\Gamma\}$ is unbounded in $\omega$.  

Let $r\leq p$ decide $g(m)$, i.e., for some $k<\omega$ suppose 
  $$r\force \dot{g}(m)=k.$$ 
Choose $\alpha\in\Gamma$ such that $f_\alpha(m)>k$, then $r$ forces a
contradiction and the Lemma is proved.
\qed

\begin{lemma}\label{lemdagg2}
Suppose $N$ is a countable standard model of ZFC+CH, $\poset$ is a countable
poset in $N$, and 
$$N\models Y\subseteq 2^\omega \mbox{ is not a  } \lambda'\mbox{ - set}$$
Then for $G$ $\poset$-generic over $N$ we have that
$$N[G]\models Y \mbox{ is not a  } \lambda'\mbox{ - set}$$
\end{lemma}
\proof
Let $D\subseteq 2^\omega$ be countable in $N$ and witness that
$Y$ is not a $\lambda'$-set, ie. there is no $G_\delta$ set $\bigcap_nU_n$
coded in $N$ with
$$\bigcap_nU_n\cap (Y\cup D)=D$$
Working in $N$ let $D=\{x_n:n<\omega\}$ and 
let $Z=Y\setminus D$ and for each $z\in Z$ define 
$f_z\in\omega^\omega$ such that $f_z(n)$ is the least $m$ such that
$x_n\res m\not= z\res m$.
Now the family $X=\{f_z:z\in Z\}$ must be unbounded in $\leq^*$ in
$N$.  Suppose not, then there exists $g\in\omega^\omega\cap N$ which
eventually dominates each element of $X$.  It follows that if we
let 
$$U_n=\bigcup_{m<n} [x_m\res n]\cup\bigcup_{m\geq n} [x_m\res g(m)]$$
then 
$$(\bigcap_{n<\omega}U_n)\cap (Y\cup D)=D$$
which is a contradiction.

It follows from Lemma \ref{lemdagg1} that $X$ is unbounded in
$N[G]$.  I claim that $D$ cannot be $G_\delta$ in $Y\cup D$ in
the model $N[G]$.  Suppose not, and let $\bigcap_{n<\omega}U_n$
be a $G_\delta$ in $N[G]$ such that 
  $$\bigcap_{n<\omega}U_n\cap(Y\cup D)=D$$
For each $n$ let $g_n\in\omega^\omega$ be such that for every $m$ we have that 
  $$[x_m\res g_n(m)]\subseteq U_n.$$
Now for any $z\in Z$ there exist a $n$ such that $z\notin U_n$.
But this means that $f_z(m)\leq g_n(m)$ for every $m$ since otherwise
  $$x_m\res g_n(m)= z\res g_n(m)$$
and then $z\in U_n$. This proves the Lemma. 
\qed

Now we prove Theorem \ref{dagger2}.  Suppose that $X\subseteq 2^\omega$ is
in $M[G]$  where $G$ is $\poset_{\omega_2}$-generic over $M$ and
 $$M[G]\models X\mbox{ is not a } \lambda' \mbox{-set}$$  
By Lowenheim-Skolem arguments there exists $\alpha<\omega_2$
such that 
$$X_\alpha=^{def}X\cap M[G_\alpha],\;\;
X_\alpha\in M[G_\alpha], \rmand 
M[G_\alpha]\models X_\alpha\mbox{ is not a } \lambda' \mbox{-set}$$
Since being a $\lambda'$-set only depends on codes for $G_\delta$-sets
and reals are added by countable suborders of 
$\poset_{[\alpha,\omega_2)}$ it follows from Lemma \ref{lemdagg2} that
 $$M[G]\models X_\alpha\mbox{ is not a } \lambda' \mbox{-set}$$
But if $f\in\omega^\omega\in M[G]$ is $\omega^{<\omega}$-generic over 
$M[G_\alpha]$ then $X_\alpha\subseteq G_f$. It follows that 
  $$M[G]\models X \mbox{ is not \daglamb}$$
as was to be proved.
\qed

\begin{flushleft}
Arnold W. Miller \\
miller@math.wisc.edu \\
http://www.math.wisc.edu/$\sim$miller\\
University of Wisconsin-Madison \\
Department of Mathematics, Van Vleck Hall \\
480 Lincoln Drive \\
Madison, Wisconsin 53706-1388 \\
\end{flushleft}

\end{document}